\documentclass[12pt]{amsart}
\usepackage{graphicx}
\usepackage{amsmath}
\usepackage{amssymb}
\usepackage{setspace}
\newtheorem{thm}{Theorem}[section]

\newtheorem{prop}[thm]{Proposition}
\theoremstyle{definition}
\newtheorem{defn}[thm]{Definition}
\theoremstyle{remark}

\theoremstyle{conclusion}

\theoremstyle{question}
\numberwithin{equation}{section}

\begin{document}
\title[The bilinear Hilbert transform and multipliers with 1D singularity]{$L^{p}$ estimates for the bilinear Hilbert transform for $1/2<p\leq2/3$: A counterexample and generalizations to non-smooth symbols}

\author{Wei Dai and Guozhen Lu}

\address{School of Mathematics and Systems Science, Beihang University (BUAA), Beijing 100191, P. R. China}
\email{daiwei@amss.ac.cn}

\address{Department of Mathematics, Wayne State University, Detroit, MI 48202, U. S. A.}
\email{gzlu@wayne.edu}

\thanks{Research of the first author  was partly supported by grants from the NNSF of China, the China Postdoctoral Science Foundation and the research of the second author was partly supported by a US NSF grant.}

\begin{abstract}
M. Lacey and C. Thiele proved in \cite{LT1} (Annals of Math. (1997)) and \cite{LT2} (Annals of Math. (1999)) that the bilinear Hilbert transform maps $L^{p_1}\times L^{p_2}\rightarrow L^{p}$ boundedly when $\frac{1}{p_1}+\frac{1}{p_2}=\frac{1}{p}$ with $1<p_{1}, \, p_{2}\leq\infty$ and $\frac{2}{3}<p<\infty$. Whether the $L^p$ estimates hold in the range $p\in (1/2,2/3]$ has remained an open problem since then. In this paper, we prove that the bilinear Hilbert transform does not map $\mathcal{F}L^{p'_{1}}\times L^{p_{2}}\rightarrow L^{p}$ for $p_1<2$ and $L^{p_{1}}\times \mathcal{F}L^{p'_{2}}\rightarrow L^{p}$ for $p_2<2$ boundedly   (Theorem \ref{main0}). In particular, this shows that the bilinear Hilbert transform neither maps $\mathcal{F}L^{p'_{1}}\times L^{p_{2}}\rightarrow L^{p}$ nor $L^{p_{1}}\times \mathcal{F}L^{p'_{2}}\rightarrow L^{p}$ for $\frac{1}{2}<p<\frac{2}{3}$. Nevertheless, we can establish $L^p$ estimates for the bilinear Fourier multipliers whose symbols are not identical to but arbitrarily close to that of the bilinear Hilbert transform in the full range $p\in(1/2,\infty)$ (Theorem \ref{main1}).
\end{abstract}
\maketitle {\small {\bf Keywords:} Bilinear Hilbert transforms; $L^{p}$ estimates; random variables; non-smooth symbols; symbols with 1-dimensional singularity; paraproducts.\\

}

\section{Introduction}

The bilinear Hilbert transform is defined by
\begin{equation}\label{1.1}
  BHT(f_{1},f_{2})(x):=p. \, v. \, \int_{\mathbb{R}}f_{1}(x-t)f_{2}(x+t)\frac{dt}{t},
\end{equation}
or equivalently, it can also be written as the bilinear multiplier operator
\begin{equation}\label{1.2}
  BHT: \,\,\, (f_{1},f_{2})\mapsto\int_{\mathbb{R}^{2}}sgn(\xi_{1}-\xi_{2})\hat{f_{1}}(\xi_{1})\hat{f_{2}}(\xi_{2})e^{2\pi ix(\xi_{1}+\xi_{2})}d\xi_{1}d\xi_{2},
\end{equation}
where $f_{1}$ and $f_{2}$ are Schwartz functions on $\mathbb{R}$.

\medskip

The boundedness problem of the bilinear Hilbert transform  was originally raised by A. P. Calder\'on in connection with the Cauchy integral along Lipschitz curves.  He conjectured that the bilinear Hilbert transform was bounded from $L^2\times L^2$ to $L^1$.
Inspired by the two classic proofs of almost everywhere convergence of Fourier series by L. Carleson \cite{Car} and C. Fefferman \cite{Fe},   and using delicate orthogonality estimates and combinatorial selection of trees, and deep time-frequency analysis, M. Lacey and C. Thiele proved  in \cite{LT1,LT2} the following celebrated $L^{p}$ estimates for the bilinear Hilbert transform.

\begin{thm}\label{BHT}(\cite{LT1,LT2})
The bilinear operator BHT maps $L^{p}(\mathbb{R})\times L^{q}(\mathbb{R})$ into $L^{r}(\mathbb{R})$ boundedly for any $1<p, \, q\leq\infty$ with $\frac{1}{p}+\frac{1}{q}=\frac{1}{r}$ and $\frac{2}{3}<r<\infty$.
\end{thm}

In general, we consider a symbol $m(\xi_{1},\xi_{2})$ that is smooth away from the singularity line $\Gamma:=\{(\xi_{1},\xi_{2})\in\mathbb{R}^{2}|\xi_{1}=\xi_{2}\}$ and satisfies
\begin{equation}\label{1.3}
  |\partial^{\alpha}m(\xi)|\lesssim\frac{1}{dist(\xi,\Gamma)^{|\alpha|}}
\end{equation}
for every $\xi=(\xi_{1},\xi_{2})\in\mathbb{R}^{2}\setminus\Gamma$ and sufficiently many multi-indices $\alpha$. Throughout this paper, $A\lesssim B$ means that there exists a universal constant $C>0$ such that $A\leq CB$. If necessary, we use explicitly $A\lesssim_{\star,\cdots,\star}B$ to indicate that there exists a positive constant $C_{\star,\cdots,\star}$ depending only on the quantities appearing in the subscript continuously such that $A\leq C_{\star,\cdots,\star}B$.

\medskip

J. Gilbert and A. Nahmod \cite{GN} proved that the $L^{p}$ estimates as $BHT$ are valid for the generalized bilinear multiplier operators $T_{m}$ associated with symbol $m$ in the same range of $p>\frac{2}{3}$ as that in \cite{LT1, LT2}.
There has been much  work related to  the bilinear operators of $BHT$ type. F. Bernicot \cite{Be} proved a pseudo-differential variant of the multiplier estimates in \cite{GN}. Uniform estimates were obtained by C. Thiele \cite{Thiele}, L. Grafakos and X. Li \cite{GL} and X. Li \cite{Li}. A two-dimensional bilinear Hilbert transform was studied by C. Demeter and C. Thiele \cite{DT}. A maximal variant of Theorem \ref{BHT} was proved by M. Lacey \cite{Lacey} and generalized by C. Demeter, T. Tao and C. Thiele \cite{DTT}. In C. Muscalu, C. Thiele and T. Tao \cite{MTT2} and J. Jung \cite{Jung}, the authors investigated various trilinear variants of the bilinear Hilbert transform. For more related results involving estimates for multi-linear singular multiplier operators, we refer to the works, e.g.,    \cite{BDNTTV, CJ, CM1, CM2, FS, GT, Jo, KS, MS, MTT1, MPTT, Thiele1} and the references therein.

\medskip

In multi-parameter cases, there is also a large amount of literature devoted to studying the estimates of multi-parameter and multi-linear operators (see \cite{CL, DL1, DL2, HL, Kesler, Luthy, LM, MS, MPTT1, MPTT2, S} and the references therein). In the bilinear and bi-parameter cases, let $\Gamma_{i}$ ($i=1,2$) be subspaces in $\mathbb{R}^{2}$, we consider operators $T_{m}$ defined by
\begin{equation}\label{1.4}
  T_{m}(f_{1},f_{2})(x):=\int_{\mathbb{R}^{4}}m(\xi,\eta)\hat{f_{1}}(\xi_{1},\eta_{1})\hat{f_{2}}(\xi_{2},\eta_{2})e^{2\pi ix\cdot((\xi_{1},\eta_{1})+(\xi_{2},\eta_{2}))}d\xi d\eta,
\end{equation}
where the symbol $m$ satisfies
\begin{equation}\label{1.5}
  |\partial_{\xi}^{\alpha}\partial_{\eta}^{\beta}m(\xi,\eta)|\lesssim\frac{1}{dist(\xi,\Gamma_{1})^{|\alpha|}}\cdot\frac{1}{dist(\eta,\Gamma_{2})^{|\beta|}}
\end{equation}
for sufficiently many multi-indices $\alpha=(\alpha_{1},\alpha_{2})$ and $\beta=(\beta_{1},\beta_{2})$. If $dim \, \Gamma_{1}=dim \, \Gamma_{2}=0$, C. Muscalu, J. Pipher, T. Tao and C. Thiele proved in \cite{MPTT1,MPTT2} that H\"{o}lder type $L^{p}$ estimates are available for $T_{m}$; J. Chen and G. Lu established in  \cite{CL} such estimates for multipliers with limited smoothness on the symbols. A multi-linear and multi-parameter pseudo-differential operator analogue has also been studied and $L^p$ estimates have been obtained by W. Dai and G. Lu  in \cite{DL1}. However, if $dim \, \Gamma_{1}=dim \, \Gamma_{2}=1$ with $\Gamma_{1}, \, \Gamma_{2}$ non-degenerate in the sense of \cite{MTT1}, let $T_{m}$ be the double bilinear Hilbert transform on polydisks $BHT\otimes BHT$ defined by
\begin{equation}\label{1.6}
  BHT\otimes BHT(f_{1},f_{2})(x,y):=p. \, v. \, \int_{\mathbb{R}^{2}}f_{1}(x-s,y-t)f_{2}(x+s,y+t)\frac{ds}{s}\frac{dt}{t},
\end{equation}
then the authors of  \cite{MPTT1}  also proved  that the operator $BHT\otimes BHT$ does not satisfy any $L^{p}$ estimates of H\"{o}lder type by constructing a counterexample. Nevertheless, under some (slightly better) logarithmic decay assumptions on the symbols, W. Dai and G. Lu proved in \cite{DL2} that the bi-parameter operators $T_{m}$ defined by \eqref{1.4}, \eqref{1.5} with singularity sets $\Gamma_{1}=\Gamma_{2}=\{(\xi_{1},\xi_{2})\in\mathbb{R}^{2}|\xi_{1}=\xi_{2}\}$ satisfy the same estimates as $BHT$. When $dim \, \Gamma_{1}=0$ and $dim \, \Gamma_{2}=1$ with $\Gamma_{2}$ non-degenerate in the sense of \cite{MTT1}, P. Silva \cite{S} and the authors of the current paper established in \cite{DL2} the $L^{p}$ estimates of H\"{o}lder type for $T_{m}$ under some conditions, which addressed the Question 8.2 in \cite{MPTT1}.

\medskip

One can observe that the $L^{p}$ estimates for the bilinear operators of $BHT$ type derived in these previous works are available only for $p>\frac{2}{3}$. In \cite{Lacey}, by constructing a counterexample, M. Lacey proved that the discrete model operators associated with the bilinear maximal functions cannot be a uniformly bounded bilinear map from $L^{p}\times L^{q}$ into $L^{r}$ for $\frac{1}{2}<r<\frac{2}{3}$. For the endpoint case $r=\frac{2}{3}$, D. Bilyk and L. Grafakos \cite{BG} proved some distributional estimates for $BHT$ of $\log$ type, then F. D. Plinio and C. Thiele \cite{PT} improved the distributional estimates by replacing the single logarithmic term with a double logarithmic term.

\medskip

Since M. Lacey and C. Thiele established the $L^p$ estimates for $\frac{2}{3}<p<\infty$ in \cite{LT1, LT2}, whether the bilinear operators of $BHT$ type satisfy $L^{p}$ estimates all the way down to $\frac{1}{2}$ has remained an open problem. Though we do not have a counterexample yet for the $L^p$ estimates for the bilinear Hilbert transform in the range of $\frac{1}{2}<p<\frac{2}{3}$, we will provide a counterexample for a modified version. We denote by $\mathcal{F}L^{p}(\mathbb{R})$ the space consisting of all functions $f$ whose Fourier transform $\hat{f}$ satisfies $\hat{f}\in L^{p}(\mathbb{R})$. The Hausdorff-Young inequality tells us that $\|\hat{f}\|_{L^{p'}(\mathbb{R})}\lesssim_{p}\|f\|_{L^{p}(\mathbb{R})}$ for $1\leq p\leq2$. Then by Theorem 1.1 due to Lacey and Thiele, it implies that the bilinear Hilbert transform maps $\mathcal{F}L^{p'_{1}}\times L^{p_{2}}\rightarrow L^{p}$ for $p_1\ge 2$  and  maps $L^{p_{1}}\times\mathcal{F}L^{p'_{2}}\rightarrow L^{p}$ for $p_2\ge 2$  with $\frac{1}{p_{1}}+\frac{1}{p_{2}}=\frac{1}{p}$. Thus  it is  interesting to ask whether the bilinear operators of $BHT$ type map $\mathcal{F}L^{p'_{1}}\times L^{p_{2}}\rightarrow L^{p}$  for $p_1<2$ or $L^{p_{1}}\times\mathcal{F}L^{p'_{2}}\rightarrow L^{p}$ for $p_2<2$ boundedly with $\frac{1}{p_{1}}+\frac{1}{p_{2}}=\frac{1}{p}$.   In this paper, we will address this issue. Our first main result in this paper is a negative answer to the boundedness of $\mathcal{F}L^{p'_{1}}\times L^{p_{2}}\rightarrow L^{p}$ for $p_1<2$  and $L^{p_{1}}\times\mathcal{F}L^{p'_{2}}\rightarrow L^{p}$ for $p_2<2$. 

\begin{thm}\label{main0}
For $p\in[\frac{1}{2},\infty)$ and $1\leq p_{1}, \, p_{2}\leq\infty$ satisfying $\frac{1}{p_{1}}+\frac{1}{p_{2}}=\frac{1}{p}$, the necessary condition for the bilinear Hilbert transform to map $\mathcal{F}L^{p'_{1}}(\mathbb{R})\times L^{p_{2}}(\mathbb{R})\rightarrow L^{p}(\mathbb{R})$ boundedly is $p_{1}\geq2$ and to map $L^{p_{1}}(\mathbb{R})\times \mathcal{F}L^{p'_{2}}(\mathbb{R})\rightarrow L^{p}(\mathbb{R})$ boundedly is $p_{2}\geq2$, respectively. In particular, the bilinear Hilbert transform maps neither $\mathcal{F}L^{p'_{1}}\times L^{p_{2}}\rightarrow L^{p}$ nor
 $L^{p_{1}}\times \mathcal{F}L^{p'_{2}}\rightarrow L^{p}$ boundedly for any $p\in[1/2,2/3)$ and $p_{1}, \, p_{2}\geq1$ satisfying $\frac{1}{p_{1}}+\frac{1}{p_{2}}=\frac{1}{p}$.
\end{thm}

Although there are no uniform $L^{p}$ estimates that are known for the bilinear Hilbert transform in the range $p\in(1/2,2/3)$, by decomposing the bilinear multiplier operator $T_{m}$ into a summation of infinitely many bilinear paraproducts without modulation invariance, we can prove that there exists a class of symbols $m$ (with one-dimensional singularity sets) which also satisfy the symbol estimates of $BHT$ type operators (see \eqref{1.3}) investigated in \cite{GN} and are arbitrarily close to the symbols of $BHT$ type operators (see \eqref{1.2} and \eqref{1.3}), such that the corresponding bilinear multiplier operators $T_{m}$ associated with symbols $m$ satisfy $L^{p}$ estimates all the way down to $\frac{1}{2}$. Our next result in this paper is the following theorem.
\begin{thm}\label{main1}
For arbitrarily given $\delta>0$, let $m(\xi_{1},\xi_{2})$ be a symbol that is smooth away from the singularity line $\Gamma:=\{(\xi_{1},\xi_{2})\in\mathbb{R}^{2}|\xi_{1}=\xi_{2}\}$ and satisfies
\begin{equation}\label{1.7}
  |\partial^{\alpha}m(\xi)|\lesssim\frac{1}{dist(\xi,\Gamma)^{|\alpha|}}\cdot
  \exp\bigg\{-\delta\big(1-\frac{|\alpha|}{3}\big)\frac{|(\xi_{1},\xi_{2})|}{dist(\xi,\Gamma)}\bigg\}, \,\,\, 0\leq|\alpha|\leq2
\end{equation}
and
\begin{equation}\label{1.8}
   |\partial^{\alpha}m(\xi)|\lesssim\frac{1}{dist(\xi,\Gamma)^{|\alpha|}}, \,\,\, |\alpha|\geq3
\end{equation}
for every $\xi=(\xi_{1},\xi_{2})\in\mathbb{R}^{2}\setminus\Gamma$ and sufficiently many multi-indices $\alpha$, then the bilinear multiplier operator $T_{m}$ defined by
\begin{equation}\label{1.9}
  T_{m}(f_{1},f_{2})(x):=\int_{\mathbb{R}^{2}}m(\xi_{1},\xi_{2})\hat{f_{1}}(\xi_{1})\hat{f_{2}}(\xi_{2})e^{2\pi ix(\xi_{1}+\xi_{2})}d\xi_{1}d\xi_{2}
\end{equation}
maps $L^{p_{1}}\times L^{p_{2}}\rightarrow L^{p}$ boundedly for any $\frac{1}{p}=\frac{1}{p_{1}}+\frac{1}{p_{2}}$ with $1<p_{1}, \, p_{2}\leq\infty$ and $\frac{1}{2}<p<\infty$. The implicit constants in the bounds depend only on $p_{1}$, $p_{2}$, $p$ when $p>\frac{2}{3}$, also depend on $\delta$ and tend to infinity as $\delta\rightarrow0$ when $\frac{1}{2}<p\leq\frac{2}{3}$.
\end{thm}

The rest of this paper is organized as follows. In Section 2 we prove Theorem \ref{main0} by using contradiction arguments. Section 3 is devoted to carrying out the proof of Theorem \ref{main1}.

\vspace{0.5cm}

{\bf Acknowledgement} The authors would like to thank Michael Lacey for his comments and for pointing out that the random variables $\epsilon_n$, given in the proof of  our original Theorem 1.2 in the first version posted in the arxiv.org,  are not independent. They also like to thank Xiaochun Li for several conversations concerning the state of affairs on the $L^p$ estimates for the bilinear Hilbert transforms throughout this work, and thank Lu Zhang, Yen Do, Loukas Grafakos and several graduate students of the second author for their comments    in the first version posted  in the arxiv.

The authors are particularly indebted to Lu Zhang for useful ideas and many helpful discussions in modifying our original counterexample  and providing many helpful ideas in constructing  the current one.  

\section{Proof of Theorem \ref{main0}}

We will prove Theorem \ref{main0} in this section by using Khinchine's inequality and contradiction arguments (see e.g., \cite{St}). To this end, without loss of generality, we first assume on the contrary that the conclusions in Theorem \ref{main0} are not true, that is, there exists some $r_{0}\in[\frac{1}{2},2)$ and $1\leq p_{0}\leq\infty$, $1\leq q_{0}<2$ such that $\frac{1}{r_{0}}=\frac{1}{p_{0}}+\frac{1}{q_{0}}$ and the bilinear Hilbert transform maps $L^{p_{0}}(\mathbb{R})\times \mathcal{F}L^{q'_{0}}(\mathbb{R})$ into $L^{r_{0}}(\mathbb{R})$ boundedly.
\begin{defn}\label{bump functions}(\cite{MS,MPTT2})
For $J\subseteq\mathbb{R}$ an arbitrary interval, we say that a smooth function $\Phi_{J}$ is a bump adapted to $J$, if and only if the following inequalities hold:
\begin{equation}\label{2.1}
  |\Phi_{J}^{(l)}(x)|\lesssim_{l,\alpha}\frac{1}{|J|^{l}}\cdot\frac{1}{\big(1+\frac{dist(x,J)}{|J|}\big)^{\alpha}}
\end{equation}
for every integer $\alpha\in\mathbb{N}$ and for sufficiently many derivatives $l\in\mathbb{N}$. If $\Phi_{J}$ is a bump adapted to $J$, we say that $|J|^{-\frac{1}{2}}\Phi_{J}$ is an $L^{2}$-normalized bump adapted to $J$.
\end{defn}
In order to get a contradiction, let us first consider two $L^{1}$-normalized even Schwartz functions $\Phi^{1}(x)$ and $\Phi^{2}(x)$ adapted to the interval $[-\frac{1}{2},\frac{1}{2}]$, such that $\hat{\Phi}^{j}\geq0$ and $supp \, \hat{\Phi}^{j}\subseteq[-\frac{1}{2},\frac{1}{2}]$ for $j=1,2$. We define Schwartz function $\Gamma(\xi,\eta):=\hat{\Phi}^{1}(\xi)\cdot\hat{\Phi}^{2}(\eta)$ for every $(\xi,\eta)\in\mathbb{R}^{2}$ such that $supp \, \Gamma\subseteq[-\frac{1}{2},\frac{1}{2}]\times[-\frac{1}{2},\frac{1}{2}]$, then define two $L^{1}$-normalized even Schwartz functions $\Psi_{1}(x)$ and $\Psi_{2}(x)$ by decomposing the Schwartz function $\Gamma(\xi,\eta)$ into a product of two $L^{\infty}$-normalized Schwartz functions $\Gamma(\xi,\eta)=:\hat{\Psi}_{1}(\xi-\eta)\cdot\hat{\Psi}_{2}\Big(\frac{\xi+\eta}{|1-|\xi-\eta||}\Big)$, such that $\hat{\Psi}_{j}\geq0$, $supp \, \hat{\Psi}_{j}\subseteq[-1,1]$ and $\hat{\Psi}_{j}\geq\frac{1}{2}$ on $[-\frac{1}{4},\frac{1}{4}]$ for $j=1,2$. In fact, the $L^{\infty}$-normalized functions $\hat{\Psi}_{2}$ may possibly be different as $\xi-\eta$ varies in $supp \, \hat{\Psi}_{1}$. However, we only need to let all these functions $\hat{\Psi}_{2}$ satisfy the property $\hat{\Psi}_{2}\geq\frac{1}{2}$ on $[-\frac{1}{4},\frac{1}{4}]$ uniformly for $\xi-\eta\in[-\frac{1}{4},\frac{1}{4}]\subseteq supp \, \hat{\Psi}_{1}$, which is enough for our proof (see \eqref{2.7} and \eqref{2.3}). Therefore, we will ignore the dependence of these functions on the variable $\xi-\eta$ hereafter and denote them by the same $L^{\infty}$-normalized Schwartz function $\hat{\Psi}_{2}$.

\medskip

Let $N$ be an arbitrarily fixed large positive integer, we define two sequences of functions
\begin{equation}\label{2.8}
   \mathcal{F}(x):=\Phi^{1}\big(x-\frac{1}{2}\big), \,\,\,\,\,\, \mathcal{G}_{N,k}(x):=\Phi^{2}\big(x-\frac{1}{2}\big)e^{4\pi i(k+2\varepsilon_{k}N)x}
\end{equation}
for every $k=1,\cdots,N$.

One can observe that
\begin{eqnarray}\label{2.4}
&&\sum_{k=1}^{N}BHT\bigg(\mathcal{F}, \, \mathcal{G}_{N,k}\bigg)(x)\\
\nonumber &=&\int_{\mathbb{R}^{2}}sgn(\xi-\eta)\sum_{k=1}^{N}e^{2\pi i(k+2\varepsilon_{k}N)}
\hat{\Phi}^{1}(\xi)\hat{\Phi}^{2}(\eta-2(k+2\varepsilon_{k}N))e^{2\pi i(x-\frac{1}{2})(\xi+\eta)}d\xi d\eta\\
\nonumber &=&-\sum_{k=1}^{N}\varepsilon_{k}\int_{\mathbb{R}^{2}}\hat{\Phi}^{1}(\xi)
\hat{\Phi}^{2}(\eta-2(k+2\varepsilon_{k}N))e^{2\pi i(x-\frac{1}{2})(\xi+\eta)}d\xi d\eta \cdot e^{2\pi i(k+2\varepsilon_{k}N)},
\end{eqnarray}
where $\{\varepsilon_{k}\}_{k=1}^{N}$ are independent and identically distributed random variables with $P(\varepsilon_{k}=\pm1)=\frac{1}{2}$ for every $k=1,\cdots,N$, i.e., $\{\varepsilon_{k}\}_{k=1}^{N}$ is a sequence with Rademacher distribution. By changing the integral variables, we can rewrite the above expression of $\sum_{k=1}^{N}BHT\bigg(\mathcal{F},\mathcal{G}_{N,k}\bigg)(x)$ in terms of Schwartz functions $\Psi_{1}$ and $\Psi_{2}$ as follows:
\begin{eqnarray}\label{2.2}
&&-\sum_{k=1}^{N}BHT\bigg(\mathcal{F}, \, \mathcal{G}_{N,k}\bigg)(x)\\
\nonumber &=&\sum_{k=1}^{N}\varepsilon_{k}\int_{\mathbb{R}^{2}}
\Gamma\bigg(\xi,\eta-2(k+2\varepsilon_{k}N)\bigg)e^{2\pi i(x-\frac{1}{2})(\xi+\eta)}d\xi d\eta \cdot e^{2\pi i(k+2\varepsilon_{k}N)} \\
\nonumber &=&\frac{1}{2}\sum_{k=1}^{N}\varepsilon_{k}\int_{\mathbb{R}^{2}}\hat{\Psi}_{1}(v)
\hat{\Psi}_{2}\Big(\frac{u}{|1-|v||}\Big)e^{2\pi i(x-\frac{1}{2})u}dudv \cdot e^{4\pi ix(k+2\varepsilon_{k}N)}\\
\nonumber &=&\frac{1}{2}\sum_{k=1}^{N}\varepsilon_{k}\int_{\mathbb{R}_{\bar{v}}}\hat{\Psi}_{1}(\bar{v})
\big|1-|\bar{v}|\big|\Psi_{2}\big(\big|1-|\bar{v}|\big|(x-\frac{1}{2})\big)d\bar{v}\cdot e^{4\pi ix(k+2\varepsilon_{k}N)}\\
\nonumber &=:&\sum_{k=1}^{N}\varepsilon_{k}\Omega(x)e^{4\pi i(k+2\varepsilon_{k}N)x},
\end{eqnarray}
where the variables $u:=\xi+\eta-2(k+2\varepsilon_{k}N)$, $v:=\xi-\eta+2(k+2\varepsilon_{k}N)$ and $\bar{v}:=\xi-\eta$ for arbitrarily large positive integer $N$.

\medskip

Now we consider arbitrary $x\in I_{N}:=(\bigcup^{\ell=-1}_{\ell=-N}I^{\ell}_{N})\cup(\bigcup^{\ell=N}_{\ell=1}I^{\ell}_{N})\subset[\frac{1}{4},\frac{3}{4}]$, where the intervals $I^{\ell}_{N}:=[\frac{1}{2}+\frac{\ell}{2^{3}N}-\frac{1}{2^{6}N},\frac{1}{2}+\frac{\ell}{2^{3}N}+\frac{1}{2^{6}N}]$ for $\ell=-N,\cdots,-1,1,\cdots,N$ and $|I^{\ell}_{N}|\simeq N^{-1}$. Since $\Psi_{1}(x)$ and $\Psi_{2}(x)$ are $L^{1}$-normalized even Schwartz functions adapted to the interval $[-\frac{1}{2},\frac{1}{2}]$ and such that $\hat{\Psi}_{j}\geq0$, $supp \, \hat{\Psi}_{j}\subseteq[-1,1]$, $\hat{\Psi}_{j}\geq\frac{1}{2}$ on $[-\frac{1}{4},\frac{1}{4}]$ for $j=1,2$, one has $\Psi_{2}(x-\frac{1}{2})\geq0$ on $I_{N}$, furthermore,
\begin{equation}\label{2.7}
\Psi_{2}(x-\frac{1}{2})=\frac{1}{2\pi}\int_{\mathbb{R}}\hat{\Psi}_{2}(\xi)\cos\big(2\pi (x-\frac{1}{2})\xi\big)d\xi
\geq\frac{1}{2\pi}\int_{-\frac{1}{4}}^{\frac{1}{4}}\frac{1}{2}\cos\frac{\pi}{4}d\xi\geq\frac{\sqrt{2}}{16\pi}
\end{equation}
for every $x\in I_{N}$, and hence we can get an estimate of lower bound for the function $\Omega$ in the right-hand side of \eqref{2.2} on $I_{N}$:
\begin{eqnarray}\label{2.3}
\big|\Omega(x)\big|&\geq&\frac{1}{2}\bigg|\int_{\mathbb{R}_{\bar{v}}}\hat{\Psi}_{1}(\bar{v})
\big|1-|\bar{v}|\big|\Psi_{2}\big(\big|1-|\bar{v}|\big|(x-\frac{1}{2})\big)d\bar{v}\bigg|\\
\nonumber &\geq& \frac{1}{4}\int^{\frac{1}{4}}_{-\frac{1}{4}}\hat{\Psi}_{1}(\bar{v})
\Psi_{2}\big(\big|1-|\bar{v}|\big|(x-\frac{1}{2})\big)d\bar{v}\geq\frac{\sqrt{2}}{256\pi}
\end{eqnarray}
for arbitrary positive integer $N$.

\medskip

Therefore, by taking average over all possible choices of i.i.d. random variables $\{\varepsilon_{k}\}_{k=1}^{N}$, we can deduce from \eqref{2.3} the following estimate for the right hand side of \eqref{2.2}:
\begin{eqnarray}\label{2.5}
  && \mathbb{E}\bigg[\bigg|\sum_{k=1}^{N}\varepsilon_{k}\Omega(x)e^{4\pi i(k+2\varepsilon_{k}N)x}\bigg|^{r_{0}}\bigg]\\
  \nonumber &\gtrsim&\big|\Omega(x)\big|^{r_{0}}\cdot\mathbb{E}\bigg[\bigg|\sum_{k=1}^{N}
  \bigg(\varepsilon_{k}\cos(4k\pi x)\cos(2^{3}N\pi x)-\sin(4k\pi x)\sin(2^{3}N\pi x)\bigg)\bigg|^{r_{0}}\bigg] \\
  \nonumber &\gtrsim& \mathbb{E}\bigg[\bigg|\sum_{k=1}^{N}\varepsilon_{k}\cos(4k\pi x)\cos(2^{3}N\pi x)\bigg|^{r_{0}}\bigg]-\bigg|\sum_{k=1}^{N}\sin(4k\pi x)\sin(2^{3}N\pi x)\bigg|^{r_{0}}
\end{eqnarray}
for every $x\in I_{N}$ and arbitrary positive integer $N$. On one hand, one can observe the fact that $|\sin(4\pi x)|\simeq(1+|\ell|)|I^{\ell}_{N}|$ and $|\cos(2^{3}N\pi x)|\simeq1$ for every $x\in I^{\ell}_{N}$ ($\ell=-N,\cdots,-1,1,\cdots,N$) and arbitrary positive integer $N$ large enough, and hence, we can deduce from Khinchine's inequality that
\begin{eqnarray}\label{2.14}
&&\mathbb{E}\bigg[\bigg|\sum_{k=1}^{N}\varepsilon_{k}\cos(4k\pi x)\cos(2^{3}N\pi x)\bigg|^{r_{0}}\bigg]\simeq\bigg(\sum_{k=1}^{N}|\cos(4k\pi x)\cos(2^{3}N\pi x)|^{2}\bigg)^{\frac{r_{0}}{2}}\\
\nonumber &\gtrsim& \bigg(N+\sum_{k=1}^{N}\cos(8k\pi x)\bigg)^{\frac{r_{0}}{2}}\gtrsim\bigg(N-\frac{1}{|\sin(4\pi x)|}\bigg)^{\frac{r_{0}}{2}}\gtrsim N^{\frac{r_{0}}{2}}-\Big(\frac{N}{1+|\ell|}\Big)^{\frac{r_{0}}{2}}
\end{eqnarray}
for every $x\in I^{\ell}_{N}$ ($\ell=-N,\cdots,-1,1,\cdots,N$) and arbitrary positive integer $N$ large enough. On the other hand, one can also observe that $|\sin(2\pi x)|\simeq(1+|\ell|)|I^{\ell}_{N}|$ for every $x\in I^{\ell}_{N}$ ($\ell=-N,\cdots,-1,1,\cdots,N$) and arbitrary positive integer $N$ large enough, and hence, we get
\begin{equation}\label{2.15}
   \bigg|\sum_{k=1}^{N}\sin(4k\pi x)\sin(2^{3}N\pi x)\bigg|^{r_{0}}\lesssim\bigg|\frac{\sin(2N\pi x)}{\sin(2\pi x)}\bigg|^{r_{0}}\lesssim\bigg(\frac{N}{1+|\ell|}\bigg)^{r_{0}},
\end{equation}
for every $x\in I^{\ell}_{N}$ ($\ell=-N,\cdots,-1,1,\cdots,N$) and arbitrary positive integer $N$ large enough.

By subtracting the estimates \eqref{2.14} and \eqref{2.15} into \eqref{2.5}, we have
\begin{equation}\label{2.16}
\mathbb{E}\bigg[\bigg|\sum_{k=1}^{N}\varepsilon_{k}\Omega(x)e^{4\pi i(k+2\varepsilon_{k}N)x}\bigg|^{r_{0}}\bigg]
   \gtrsim N^{\frac{r_{0}}{2}}-\Big(\frac{N}{1+|\ell|}\Big)^{r_{0}}
\end{equation}
for every $x\in I^{\ell}_{N}$ ($\ell=-N,\cdots,-1,1,\cdots,N$) and arbitrary positive integer $N$ large enough.

\medskip

Now we define function
\begin{equation}\label{2.12}
  \mathcal{G}_{N}(x):=\sum_{k=1}^{N}\mathcal{G}_{N,k}(x)
\end{equation}
for arbitrarily large positive integer $N$. Suppose that the bilinear Hilbert transform maps $L^{p_{0}}(\mathbb{R})\times\mathcal{F}L^{q'_{0}}(\mathbb{R})\rightarrow L^{r_{0}}(\mathbb{R})$ boundedly, where $r_{0}\in[\frac{1}{2},2)$ and $1\leq p_{0}\leq\infty$, $1\leq q_{0}<2$ such that $\frac{1}{r_{0}}=\frac{1}{p_{0}}+\frac{1}{q_{0}}$, so we can obtain that
\begin{equation}\label{2.0}
\bigg\|\sum_{k=1}^{N}BHT\bigg(\mathcal{F}, \, \mathcal{G}_{N,k}\bigg)\bigg\|^{r_{0}}_{L^{r_{0}}(I_{N})}\lesssim
\big\|\mathcal{F}\big\|^{r_{0}}_{L^{p_{0}}}\cdot\big\|\widehat{\mathcal{G}_{N}}\big\|^{r_{0}}_{L^{q'_{0}}}
\end{equation}
for arbitrarily large positive integer $N$.

\medskip

One easily observes that $\|\mathcal{F}\|_{L^{p_{0}}}\lesssim1$, and also has the following estimate for $\mathcal{G}_{N}$:
\begin{eqnarray}\label{2.10}
  \|\widehat{\mathcal{G}_{N}}\|^{q'_{0}}_{L^{q'_{0}}}&\lesssim&\int_{\mathbb{R_{\xi}}}
  \bigg|\mathcal{F}\Big(\Phi^{2}(x-\frac{1}{2})\sum_{k=1}^{N}e^{4\pi i(k+2\varepsilon_{k}N)x}\Big)(\xi)\bigg|^{q'_{0}}d\xi \\
 \nonumber &\lesssim& \int_{\mathbb{R}_{\xi}}\bigg|\sum_{k=1}^{N}\hat{\Phi}^{2}(\xi-2k-2^{2}\varepsilon_{k}N)e^{-\pi i(\xi-2(k+2\varepsilon_{k}N))}\bigg|^{q'_{0}}d\xi \\
 \nonumber &\lesssim& \sum_{k=1}^{N}\int_{\mathbb{R}_{\xi}}\big|\hat{\Phi}^{2}(\xi)\big|^{q'_{0}}d\xi\lesssim N
\end{eqnarray}
for arbitrarily large positive integer $N$.

\medskip

Then we can deduce from \eqref{2.0} and \eqref{2.10} that the following upper bounds
\begin{equation}\label{2.6}
  \mathbb{E}\bigg[\bigg\|\sum_{k=1}^{N}BHT\bigg(\mathcal{F}, \, \mathcal{G}_{N,k}\bigg)\bigg\|^{r_{0}}_{L^{r_{0}}(I_{N})}\bigg]\lesssim N^{(1-\frac{1}{q_{0}})r_{0}}
\end{equation}
hold true with bounds that are uniform with respect to arbitrary large positive integer $N$, while the estimates \eqref{2.2} and \eqref{2.16} yield that
\begin{eqnarray}\label{2.9}
&& \mathbb{E}\bigg[\bigg\|\sum_{k=1}^{N}BHT\bigg(\mathcal{F}, \, \mathcal{G}_{N,k}\bigg)\bigg\|^{r_{0}}_{L^{r_{0}}(I_{N})}\bigg]
=\int_{I_{N}}\mathbb{E}\bigg[\bigg|\sum_{k=1}^{N}\varepsilon_{k}\Omega(x)e^{4\pi i(k+2\varepsilon_{k}N)x}\bigg|^{r_{0}}\bigg]dx\\
\nonumber &\gtrsim& \sum_{\ell=-N, \, \ell\neq0}^{\ell=N}\bigg(N^{\frac{r_{0}}{2}}-\Big(\frac{N}{1+|\ell|}\Big)^{r_{0}}\bigg)\cdot\big|I^{\ell}_{N}\big|
\gtrsim N^{\frac{r_{0}}{2}}-\sum_{|\ell|=1}^{N}\frac{N^{r_{0}-1}}{(1+|\ell|)^{r_{0}}}\\
\nonumber &\gtrsim& N^{\frac{r_{0}}{2}}-N^{r_{0}-1}(1+N)^{1-r_{0}}\gtrsim N^{\frac{r_{0}}{2}}
\end{eqnarray}
for arbitrary positive integer $N$ large enough. Combining the estimates \eqref{2.6} and \eqref{2.9}, we must have
\begin{equation}\label{2.11}
  N^{\frac{r_{0}}{2}}\lesssim N^{(1-\frac{1}{q_{0}})r_{0}}
\end{equation}
for arbitrary positive integer $N$ large enough, which implies that $\frac{r_{0}}{2}\leq(1-\frac{1}{q_{0}})r_{0}$ holds true, that is, we must have $q_{0}\geq2$. Thus, the necessary condition for the bilinear Hilbert transform to map $L^{p}\times\mathcal{F}L^{q'}\rightarrow L^{r}$ boundedly is $q\geq2$. For any $r\in[1/2,2/3)$ and $p, \, q\geq1$ satisfying $\frac{1}{p}+\frac{1}{q}=\frac{1}{r}$, one can easily observe that both $p, \, q\in[1,2)$, therefore the bilinear Hilbert transform does not map $\mathcal{F}L^{p'}\times L^{q}\rightarrow L^{r}$ and $L^{p}\times \mathcal{F}L^{q'}\rightarrow L^{r}$ boundedly. This concludes our proof of Theorem \ref{main0}.

\section{Proof of Theorem \ref{main1}}

For arbitrarily given $\delta>0$, one can observe that the symbols $m$ defined by \eqref{1.7} and \eqref{1.8} also satisfy the estimates
\begin{equation}\label{3.1}
  |\partial^{\alpha}m(\xi)|\lesssim\frac{1}{|\xi_{1}-\xi_{2}|^{|\alpha|}}
\end{equation}
for every $\xi=(\xi_{1},\xi_{2})\in\mathbb{R}^{2}\setminus\Gamma$ and sufficiently many multi-indices $\alpha$. Therefore, we deduce from \cite{GN} that the bilinear operators $T_{m}$ given by \eqref{1.9} satisfy $L^{p}$ estimates of H\"{o}lder type for $p>\frac{2}{3}$, and the implicit constants in the bounds depend only on $p_{1}$, $p_{2}$, $p$. In this section, we will focus on proving $L^{p}$ estimates of $T_{m}$ for $p\leq\frac{2}{3}$.

\subsection{Decomposition into a summation of infinitely many bilinear multipliers}
As we can see from the study of multi-parameter and multi-linear Coifman-Meyer multiplier operators (see e.g. \cite{MTT1,MPTT1,MTT2,MPTT2}), a standard approach to obtain $L^{p}$ estimates of the  bilinear operators $T_{m}$ is to reduce it into discrete sums of inner products with wave packets (see \cite{Thiele1}).

\medskip

First, we need to decompose the symbol $m(\xi)$ in a natural way. To this end, we will decompose the region $\{\xi=(\xi_{1},\xi_{2})\in\mathbb{R}^{2}: \xi_{1}\neq\xi_{2}\}$ by using \emph{Whitney squares} with respect to the singularity point $\{\xi_{1}=\xi_{2}=0\}$. In order to describe our discretization procedure clearly, let us first recall some standard notation and definitions in \cite{MTT2}.

\medskip

An interval $I$ on the real line $\mathbb{R}$ is called dyadic if it is of the form $I=2^{-k}[n, \, n+1]$ for some $k, \, n\in\mathbb{Z}$. An interval is said to be a \emph{shifted dyadic interval} if it is of the form $2^{-k}[j+\alpha,j+1+\alpha]$ for any $k,j\in\mathbb{Z}$ and $\alpha\in\{0,\frac{1}{3},-\frac{1}{3}\}$. A \emph{shifted dyadic cube} is a set of the form $Q=Q_{1}\times Q_{2}\times Q_{3}$, where each $Q_{j}$ is a shifted dyadic interval and they all have the same length. A \emph{shifted dyadic quasi-cube} is a set $Q=Q_{1}\times Q_{2}\times Q_{3}$, where $Q_{j}$ ($j=1,2,3$) are shifted dyadic intervals satisfying less restrictive condition $|Q_{1}|\simeq|Q_{2}|\simeq|Q_{3}|$. One easily observes that for every cube $Q\subseteq\mathbb{R}^{3}$, there exists a shifted dyadic cube $\widetilde{Q}$ such that $Q\subset\frac{7}{10}\widetilde{Q}$ (the cube having the same center as $\widetilde{Q}$ but with side length $\frac{7}{10}$ that of $\widetilde{Q}$) and $diam(Q)\simeq diam(\widetilde{Q})$.

\medskip

The same terminology will also be used in the plane $\mathbb{R}^{2}$. The only difference is that the previous cubes now become squares.

For any cube and square $Q$, we will denote the side length of $Q$ by $\ell(Q)$ for short and denote the reflection of $Q$ with respect to the origin by $-Q$ hereafter.

\medskip

By writing the characteristic function $1_{\xi_{1}\neq\xi_{2}}$ of the region $\{\xi\in\mathbb{R}^{2}: \, \xi_{1}\neq\xi_{2}\}$ into finite sum of smoothed versions of characteristic functions of the cones $\{\xi_{2}>|\xi_{1}|\}$, $\{\xi_{2}<-|\xi_{1}|\}$, $\{\xi_{1}>|\xi_{2}|\}$ and $\{\xi_{1}<-|\xi_{2}|\}$, we can decompose the bilinear operator $T_{m}$ into a finite sum of four parts. Since all the operators obtained in this decomposition can be treated in the same way, without loss of generalizations, we will discuss in detail only one of them hereafter, for instance, the bilinear operator $T_{m_{\mathcal{A}}}$ given by smoothly truncating the symbol $m$ on the cone $\mathcal{A}:=\{\xi=(\xi_{1},\xi_{2})\in\mathbb{R}^{2}: \, \xi_{2}>|\xi_{1}|\}$.

\medskip

For this purpose, we consider the collection $\mathcal{Q}_{n}$ of all shifted dyadic squares $Q=Q_{1}\times Q_{2}$ satisfying the property that
\begin{equation}\label{3.2}
  Q\cap\Gamma=\emptyset \,\,\,\,\,\,\,\, \text{and} \,\,\,\,\,\,\,\, diam(Q)\simeq10^{-4}\bigg(\frac{1}{\sqrt{n+1}}-\frac{1}{\sqrt{n+2}}\bigg)\cdot dist(Q,(0,0))
\end{equation}
for arbitrary positive integer $n=1,2,\cdots$, where the singularity line $\Gamma=\{\xi=(\xi_{1},\xi_{2})\in\mathbb{R}^{2}: \, \xi_{1}=\xi_{2}\}$. Correspondingly, we decompose the cone $\mathcal{A}$ into a sequence of cones $\mathcal{C}_{1},\mathcal{C}_{2},\cdots,\mathcal{C}_{n},\cdots$ defined by
\begin{equation}\label{3.3}
  \mathcal{C}_{1}:=\{\xi=(\xi_{1},\xi_{2})\in\mathbb{R}^{2}: \, \frac{\sqrt{3}}{2}\leq\frac{dist(\xi,\Gamma)}{|(\xi_{1},\xi_{2})|}\leq1, \, \xi_{2}>|\xi_{1}|\},
\end{equation}
\begin{equation}\label{3.4}
  \mathcal{C}_{2}:=\{\xi=(\xi_{1},\xi_{2})\in\mathbb{R}^{2}: \, \frac{1}{2}\leq\frac{dist(\xi,\Gamma)}{|(\xi_{1},\xi_{2})|}\leq\frac{\sqrt{3}}{2}, \, \xi_{2}>|\xi_{1}|\}
\end{equation}
and
\begin{equation}\label{3.5}
  \mathcal{C}_{n}:=\{\xi=(\xi_{1},\xi_{2})\in\mathbb{R}^{2}: \, \frac{1}{\sqrt{n+2}}\leq\frac{dist(\xi,\Gamma)}{|(\xi_{1},\xi_{2})|}
  \leq\frac{1}{\sqrt{n+1}}, \, \xi_{2}>|\xi_{1}|\}
\end{equation}
for every $n\geq3$. The borderlines $\{\mathcal{L}_{n}\}_{n=1}^{\infty}$ of these cones are defined by
\begin{equation}\label{3.6}
  \mathcal{L}_{1}:=\{\xi=(\xi_{1},\xi_{2})\in\mathbb{R}^{2}: \, \xi_{2}=\bigg(\tan\frac{7}{12}\pi\bigg)\xi_{1}>|\xi_{1}|\},
\end{equation}
\begin{equation}\label{3.7}
  \mathcal{L}_{2}:=\{\xi=(\xi_{1},\xi_{2})\in\mathbb{R}^{2}: \, \xi_{2}=\bigg(\tan\frac{5}{12}\pi\bigg)\xi_{1}>|\xi_{1}|\}
\end{equation}
and
\begin{equation}\label{3.8}
  \mathcal{L}_{n}:=\{\xi=(\xi_{1},\xi_{2})\in\mathbb{R}^{2}: \, \xi_{2}=\bigg(\tan\big(\frac{\pi}{4}+\arcsin\frac{1}{\sqrt{n+2}}\big)\bigg)\xi_{1}>|\xi_{1}|\}
\end{equation}
for every $n\geq3$.

\medskip

Now let us define the following disjoint collections of shifted dyadic squares
\begin{equation}\label{3.9}
  \mathbb{Q}_{n}:=\{Q=Q_{1}\times Q_{2}: \, Q\in\mathcal{Q}_{n}, \,\,\,  Q\cap\mathcal{C}_{n}\neq\emptyset \,\,\, \text{and} \,\,\, Q\cap\mathcal{L}_{n}=\emptyset\}
\end{equation}
for $n=1,2,\cdots$. Since the set of squares $\{\frac{7}{10}Q \, | \, Q\in\bigcup_{n=1}^{\infty}\mathbb{Q}_{n}\}$ also forms a finitely overlapping cover of the cone $\mathcal{A}:=\bigcup_{n=1}^{\infty}\mathcal{C}_{n}=\{\xi=(\xi_{1},\xi_{2}): \, |\xi_{1}|<\xi_{2}\}$, by a standard partition of unity, we can write the smoothed characteristic function $\tilde{\chi}_{\{|\xi_{1}|<\xi_{2}\}}$ of the cone $\mathcal{A}$ as
\begin{equation}\label{3.10}
  \tilde{\chi}_{\{|\xi_{1}|<\xi_{2}\}}=\sum_{n=1}^{\infty}\sum_{Q\in\mathbb{Q}_{n}}\phi_{Q}(\xi_{1},\xi_{2}),
\end{equation}
where each $\phi_{Q}$ is a smooth bump function adapted to $Q$ and supported in $\frac{8}{10}Q$. One observes that all the collections $\mathbb{Q}_{n}$ ($n=1,2,\cdots$) can be decomposed further into at most $10^{8}n^{\frac{3}{2}}$ (modulo some fixed constant $C_{0}$ that is independent of $n=1,2,\cdots$) disjoint sub-collections $\mathbb{Q}^{k}_{n}$ ($1\leq k\leq C_{0}10^{8}n^{\frac{3}{2}}$) which contains only one unique shifted dyadic square $Q$ with the fixed scale $\ell(Q)=2^{l_{0}}$ for some arbitrarily given integer $l_{0}\in\mathbb{Z}$, that is,
\begin{equation}\label{3.11}
  \mathbb{Q}_{n}=\bigcup_{k=1}^{\sim10^{8}n^{\frac{3}{2}}}\mathbb{Q}^{k}_{n}
\end{equation}
for every $n=1,2,\cdots$.

\medskip

In order to prove $L^{p}$ estimates for $T_{m}$ (Theorem \ref{main1}), it's enough for us to investigate the bilinear operator $T_{m_{\mathcal{A}}}$ given by
\begin{equation}\label{3.12}
  T_{m_{\mathcal{A}}}(f_{1},f_{2})(x):=\int_{\mathbb{R}^{2}}m(\xi_{1},\xi_{2})\tilde{\chi}_{\{|\xi_{1}|<\xi_{2}\}}(\xi_{1},\xi_{2})
  \hat{f_{1}}(\xi_{1})\hat{f_{2}}(\xi_{2})e^{2\pi ix(\xi_{1}+\xi_{2})}d\xi_{1}d\xi_{2}.
\end{equation}
By using \eqref{3.10} and \eqref{3.11}, we can decompose $T_{m_{\mathcal{A}}}$ into a summation of bilinear multiplier operators:
\begin{equation}\label{3.13}
  T_{m_{\mathcal{A}}}=\sum_{n=1}^{N}\sum_{k=1}^{\sim10^{8}n^{\frac{3}{2}}}T^{k}_{m_{\mathcal{A}},n},
\end{equation}
where the bilinear operators $T^{k}_{m_{\mathcal{A}},n}$ are given by
\begin{equation}\label{3.14}
  T^{k}_{m_{\mathcal{A}},n}(f_{1},f_{2})(x):=\sum_{Q\in\mathbb{Q}^{k}_{n}}\int_{\mathbb{R}^{2}}m(\xi_{1},\xi_{2})\phi_{Q}(\xi_{1},\xi_{2})
  \hat{f_{1}}(\xi_{1})\hat{f_{2}}(\xi_{2})e^{2\pi ix(\xi_{1}+\xi_{2})}d\xi_{1}d\xi_{2}
\end{equation}
for every $n=1,2,\cdots$ and $1\leq k\lesssim10^{8}n^{\frac{3}{2}}$. Therefore, the proof of Theorem \ref{main1} can be reduced to proving $L^{p}$ estimates for each single bilinear multipliers $T^{k}_{m_{\mathcal{A}},n}$ with the constants in the bounds being independent of $k$ and having enough decay which is acceptable for summation with respect to $n$, that is, the following proposition.
\begin{prop}\label{single multiplier}
For every $\delta>0$, $n=1,2,\cdots$ and $1\leq k\lesssim10^{8}n^{\frac{3}{2}}$, let $T^{k}_{m_{\mathcal{A}},n}$ be the bilinear multiplier operator defined by \eqref{3.14} with the symbol $m$ satisfies the differential estimates \eqref{1.7} and \eqref{1.8}, then we have
\begin{equation}\label{3.15}
   \|T^{k}_{m_{\mathcal{A}},n}(f_{1},f_{2})\|_{L^{p}(\mathbb{R})}\lesssim_{n,p,p_{1},p_{2},\delta}\|f_{1}\|_{L^{p_{1}}(\mathbb{R})}
   \cdot\|f_{2}\|_{L^{p_{2}}(\mathbb{R})},
\end{equation}
provided that $\frac{1}{p}=\frac{1}{p_{1}}+\frac{1}{p_{2}}$ with $1<p_{1}, \, p_{2}\leq\infty$ and $\frac{1}{2}<p<\infty$. Moreover, the implicit constants in the bounds satisfy
\begin{equation}\label{3.16}
  C_{n,p,p_{1},p_{2},\delta}\lesssim_{p,p_{1},p_{2}}e^{-\delta\sqrt{n}}+n^{-3}.
\end{equation}
\end{prop}

\subsection{Reduce each single bilinear multiplier to a discrete model operator}

Now we consider arbitrarily fixed $n=1,2,\cdots$ and $1\leq k\lesssim10^{8}n^{\frac{3}{2}}$. For each shifted dyadic square $Q\in\mathbb{Q}^{k}_{n}$, one observes that there exist bump functions $\phi_{Q_{i},i}$ ($i=1,2$) adapted to the shifted dyadic interval $Q_{i}$ such that $supp \, \phi_{Q_{i},i}\subseteq\frac{9}{10}Q_{i}$ and $\phi_{Q_{i},i}\equiv1$ on $\frac{8}{10}Q_{i}$ ($i=1,2$) respectively. Notice that $supp \, \phi_{Q}\subseteq\frac{8}{10}Q$, thus one has $\phi_{Q_{1},1}\cdot\phi_{Q_{2},2}\equiv1$ on $supp \, \phi_{Q}$. Since $\xi_{1}\in supp \, \phi_{Q_{1},1}\subseteq\frac{9}{10}Q_{1}$ and $\xi_{2}\in supp \, \phi_{Q_{2},2}\subseteq\frac{9}{10}Q_{2}$, it follows that $-\xi_{1}-\xi_{2}\in-\frac{9}{10}Q_{1}-\frac{9}{10}Q_{2}$, and as a consequence, one can find a shifted dyadic interval $Q_{3}$ with the property that $-\frac{9}{10}Q_{1}-\frac{9}{10}Q_{2}\subseteq\frac{7}{10}Q_{3}$ and satisfying $|Q_{1}|=|Q_{2}|\simeq|Q_{3}|$. In particular, there exists bump function $\phi_{Q_{3},3}$ adapted to $Q_{3}$ and supported in $\frac{9}{10}Q_{3}$ such that $\phi_{Q_{3},3}\equiv1$ on $-\frac{9}{10}Q_{1}-\frac{9}{10}Q_{2}$.

\medskip

We denote by $\mathbf{Q}^{k}_{n}$ the collection of all shifted dyadic quasi-cubes $Q:=Q_{1}\times Q_{2}\times Q_{3}$ with $Q_{1}\times Q_{2}\in\mathbb{Q}^{k}_{n}$ and $Q_{3}$ be defined as above. Assuming this we then observe that, for any $Q$ in such a collection $\mathbf{Q}^{k}_{n}$, there exists a unique shifted dyadic cube $\widetilde{Q}$ in $\mathbb{R}^{3}$ such that $Q\subseteq\frac{7}{10}\widetilde{Q}$ and with property that $diam(Q)\simeq diam(\widetilde{Q})$. This allows us in particular to assume further that $\mathbf{Q}^{k}_{n}$ is a collection of shifted dyadic cubes (that is, $|Q_{1}|=|Q_{2}|=|Q_{3}|=\ell(Q)$).

\medskip

Now consider the trilinear form $\Lambda^{k}_{m_{\mathcal{A}},n}(f_{1},f_{2},f_{3})$ associated to $T^{k}_{m_{\mathcal{A}},n}(f_{1},f_{2})$, which can be written as
\begin{eqnarray}\label{3.17}
&&\Lambda^{k}_{m_{\mathcal{A}},n}(f_{1},f_{2},f_{3}):=\int_{\mathbb{R}^{2}}T^{k}_{m_{\mathcal{A}},n}(f_{1},f_{2})(x)f_{3}(x)dx\\
\nonumber &=&\sum_{Q\in\mathbf{Q}^{k}_{n}}\int_{\xi_{1}+\xi_{2}+\xi_{3}=0}
 m_{Q}(\xi_{1},\xi_{2},\xi_{3})(f_{1}\ast\check{\phi}_{Q_{1},1})^{\wedge}(\xi_{1})(f_{2}\ast\check{\phi}_{Q_{2},2})^{\wedge}(\xi_{2})
 (f_{3}\ast\check{\phi}_{Q_{3},3})^{\wedge}(\xi_{3})d\xi,
\end{eqnarray}
where $\xi=(\xi_{1},\xi_{2},\xi_{3})\in\mathbb{R}^{3}$, while
\begin{equation}\label{3.18}
  m_{Q}(\xi_{1},\xi_{2},\xi_{3}):=m(\xi_{1},\xi_{2})\cdot\phi_{Q_{1}\times Q_{2}}(\xi_{1},\xi_{2})\cdot\widetilde{\phi}_{Q_{3},3}(\xi_{3}),
\end{equation}
where the function $\phi_{Q_{1}\times Q_{2}}(\xi_{1},\xi_{2})$ is one term of the partition of unity defined in \eqref{3.10}, $\widetilde{\phi}_{Q_{3},3}$ is an appropriate smooth function of variable $\xi_{3}$ supported on a slightly larger interval (with a constant magnification independent of $\ell(Q)$) than $supp \, \phi_{Q_{3},3}$, which equals $1$ on $supp \, \phi_{Q_{3},3}$. We can decompose $m_{Q}(\xi_{1},\xi_{2},\xi_{3})$ as a Fourier series:
\begin{equation}\label{3.19}
  m_{Q}(\xi_{1},\xi_{2},\xi_{3})=\sum_{n_{1},n_{2},n_{3}\in\mathbb{Z}}C^{\delta, \, Q}_{n_{1},n_{2},n_{3}}
  e^{2\pi i(n_{1},n_{2},n_{3})\cdot(\xi_{1},\xi_{2},\xi_{3})/\ell(Q)},
\end{equation}
where the Fourier coefficients $C^{\delta, \, Q}_{n_{1},n_{2},n_{3}}$ are given by
\begin{equation}\label{3.20}
C^{\delta, \, Q}_{n_{1},n_{2},n_{3}}=\int_{\mathbb{R}^{3}}m_{Q}(\ell(Q)\xi_{1},\ell(Q)\xi_{2},\ell(Q)\xi_{3})
e^{-2\pi i(n_{1}\xi_{1}+n_{2}\xi_{2}+n_{3}\xi_{3})}d\xi_{1}d\xi_{2}d\xi_{3}
\end{equation}
for every shifted dyadic cube $Q\in\mathbf{Q}^{k}_{n}$. Then, by a straightforward calculation, we can rewrite \eqref{3.17} as
\begin{eqnarray}\label{3.21}
\nonumber &&\Lambda^{k}_{m_{\mathcal{A}},n}(f_{1},f_{2},f_{3})=\sum_{Q\in\mathbf{Q}^{k}_{n}}
\sum_{n_{1},n_{2},n_{3}\in\mathbb{Z}}C^{\delta, \, Q}_{n_{1},n_{2},n_{3}}\int_{\mathbb{R}}(f_{1}\ast\check{\phi}_{Q_{1},1})(x-\frac{n_{1}}{\ell(Q)})\\
&&\quad\quad\quad\quad\quad\quad\quad\quad\quad\quad\quad\quad\quad\quad
\times(f_{2}\ast\check{\phi}_{Q_{2},2})(x-\frac{n_{2}}{\ell(Q)})(f_{3}\ast\check{\phi}_{Q_{3},3})(x-\frac{n_{3}}{\ell(Q)})dx.
\end{eqnarray}

\begin{defn}\label{wave packet}(\cite{MTT2,Thiele1})
An arbitrary dyadic rectangle of area $1$ in the phase-space plane is called a \emph{Heisenberg box} or \emph{tile}. Let $P:=I_{P}\times\omega_{P}$ be a tile. A $L^{2}$-normalized wave packet on $P$ is a function $\Phi_{p}$ which has Fourier support $supp \, \hat{\Phi}_{p}\subseteq\frac{9}{10}\omega_{P}$ and obeys the estimates
\begin{equation*}
  |\Phi_{P}(x)|\lesssim|I_{P}|^{-\frac{1}{2}}\bigg(1+\frac{dist(x,I_{P})}{|I_{P}|}\bigg)^{-M}
\end{equation*}
for all $M>0$, where the implicit constant depends on $M$.
\end{defn}

Now we define $\phi^{n_{i}}_{Q_{i},i}(\xi_{i}):=e^{2\pi in_{i}\xi_{i}/\ell(Q)}\cdot\phi_{Q_{i},i}(\xi_{i})$ for $i=1,2,3$. By the construction of the collection $\mathbf{Q}^{k}_{n}$ of shifted dyadic cubes, there exists only one unique cube $Q\in\mathbf{Q}^{k}_{n}$ such that $\ell(Q)=|Q_{1}|=|Q_{2}|=|Q_{3}|=2^{l}$ for every $l\in\mathbb{Z}$. By splitting the real line $\mathbb{R}$ into disjoint union of unit intervals, performing the $L^{2}$-normalization procedure and simple calculations, we can rewrite \eqref{3.21} as
\begin{eqnarray}\label{3.22}
   &&\Lambda^{k}_{m_{\mathcal{A}},n}(f_{1},f_{2},f_{3})\\
  \nonumber &=&\sum_{n_{1},n_{2},n_{3}\in\mathbb{Z}}\sum_{Q\in\mathbf{Q}^{k}_{n}}
  \int_{0}^{1}\sum_{\substack{I \,\, \text{dyadic}, \\ |I|=\ell(Q)^{-1}=2^{-l}}}\frac{C^{\delta,Q}_{n_{1},n_{2},n_{3}}}{|I|^{\frac{1}{2}}}
  \langle f_{1},\check{\phi}^{n_{1},\nu}_{I,Q_{1},1}\rangle\langle f_{2},\check{\phi}^{n_{2},\nu}_{I,Q_{2},2}\rangle
  \langle f_{3},\check{\phi}^{n_{3},\nu}_{I,Q_{3},3}\rangle d\nu\\
 \nonumber &=:&\sum_{n_{1},n_{2},n_{3}\in\mathbb{Z}}\int_{0}^{1}\sum_{P\in\mathbb{P}^{k}_{n}}
 \frac{C^{\delta,Q_{P}}_{n_{1},n_{2},n_{3}}}{|I_{P}|^{\frac{1}{2}}}
 \langle f_{1},\Phi^{1,n_{1},\nu}_{P_{1}}\rangle\langle f_{2},\Phi^{2,n_{2},\nu}_{P_{2}}\rangle
 \langle f_{3},\Phi^{3,n_{3},\nu}_{P_{3}}\rangle d\nu,
\end{eqnarray}
where the notation $\langle\cdot,\cdot\rangle$ denotes the complex scalar $L^{2}$ inner product, the Fourier coefficients $C^{\delta,Q_{P}}_{n_{1},n_{2},n_{3}}:=C^{\delta,Q}_{n_{1},n_{2},n_{3}}$, the \emph{tri-tiles} $P:=(P_{1},P_{2},P_{3})$, the \emph{tiles} $P_{j}:=I_{P_{j}}\times\omega_{P_{j}}$ with \emph{time intervals} $I_{P_{j}}:=I=2^{-l}[m, m+1]=:I_{P}$ and the \emph{frequency intervals} $\omega_{P_{j}}:=Q_{j}$ for $j=1,2,3$, the frequency cubes $Q_{P}:=Q=\omega_{P_{1}}\times\omega_{P_{2}}\times\omega_{P_{3}}\in\mathbf{Q}^{k}_{n}$, $\mathbb{P}^{k}_{n}$ denotes the collection of such tri-tiles $P$ with frequency cubes $Q_{P}\in\mathbf{Q}^{k}_{n}$, while the $L^{2}$-normalized wave packets $\Phi^{i,n_{i},\nu}_{P_{i}}$ associated with the Heisenberg boxes $P_{i}$ are defined by $\Phi^{i,n_{i},\nu}_{P_{i}}(x):=\check{\phi}^{n_{i},\nu}_{I,Q_{i},i}(x):=2^{-\frac{l}{2}}\overline{\check{\phi}^{n_{i}}_{Q_{i},i}(2^{-l}(m+\nu)-x)}$ for $i=1,2,3$.

\medskip

By taking advantage of the differential estimates \eqref{1.7} and \eqref{1.8} for symbol $m(\xi_{1},\xi_{2})$, one deduces from the expression of Fourier coefficients \eqref{3.20} and integrating by parts sufficiently many times that
\begin{eqnarray}\label{FC}
 \nonumber |C^{\delta,Q_{P}}_{n_{1},n_{2},n_{3}}|&\lesssim&\prod_{j=1}^{3}\frac{1}{(1+|n_{j}|)^{1000}}\cdot\{e^{-\delta\sqrt{n+1}}
  +\frac{e^{-\frac{2\delta}{3}\sqrt{n+1}}}{n+1}+\frac{e^{-\frac{\delta}{3}\sqrt{n+1}}}{(n+1)^{2}}+\sum_{i=3}^{3000}\frac{1}{(n+1)^{i}}\}\\
  &\lesssim& \prod_{j=1}^{3}\frac{1}{(1+|n_{j}|)^{1000}}\cdot\{e^{-\delta\sqrt{n}}+n^{-3}\}
\end{eqnarray}
for any tri-tiles $P\in\mathbb{P}^{k}_{n}$.

\medskip

Observe that the rapid decay with respect to the parameters $n_{1},n_{2},n_{3}\in\mathbb{Z}$ in \eqref{FC} is acceptable for summation, all the functions $\Phi^{j,n_{j},\nu}_{P_{j}}$ ($j=1,2,3$) are $L^{2}$-normalized and are wave packets associated with the Heisenberg boxes $P_{j}$ uniformly with respect to the parameters $n_{j}$, therefore we only need to consider from now on the part of the trilinear form $\Lambda^{k}_{m_{\mathcal{A}},n}(f_{1},f_{2},f_{3})$ defined in \eqref{3.22} corresponding to $n_{1}=n_{2}=n_{3}=0$:
\begin{equation}\label{3.23}
  \dot{\Lambda}^{k}_{m_{\mathcal{A}},n}(f_{1},f_{2},f_{3}):=\int_{0}^{1}\sum_{P\in\mathbb{P}^{k}_{n}}
 \frac{C^{\delta}_{Q_{P}}}{|I_{P}|^{\frac{1}{2}}}
 \langle f_{1},\Phi^{1,\nu}_{P_{1}}\rangle\langle f_{2},\Phi^{2,\nu}_{P_{2}}\rangle
 \langle f_{3},\Phi^{3,\nu}_{P_{3}}\rangle d\nu,
\end{equation}
where $C^{\delta}_{Q_{P}}:=C^{\delta,Q_{P}}_{0,0,0}$ and $\Phi^{i,\nu}_{P_{i}}:=\Phi^{i,0,\nu}_{P_{i}}$ for $\nu\in[0,1]$ and $i=1,2,3$.

\medskip

The bilinear operator corresponding to the trilinear form $\dot{\Lambda}^{k}_{m_{\mathcal{A}},n}(f_{1},f_{2},f_{3})$ can be written as
\begin{equation}\label{3.24}
  \dot{\Pi}^{\delta}_{\mathbb{P}^{k}_{n}}(f_{1},f_{2})(x)=\int_{0}^{1}\sum_{P\in\mathbb{P}^{k}_{n}}
 \frac{C^{\delta}_{Q_{P}}}{|I_{P}|^{\frac{1}{2}}}
 \langle f_{1},\Phi^{1,\nu}_{P_{1}}\rangle\langle f_{2},\Phi^{2,\nu}_{P_{2}}\rangle\Phi^{3,\nu}_{P_{3}}(x)d\nu.
\end{equation}

Since $\dot{\Pi}^{\delta}_{\mathbb{P}^{k}_{n}}(f_{1},f_{2})$ is an average of some discrete bilinear model operators depending on the parameters $\nu\in[0,1]$, it is enough to prove the $L^{p}$ estimates of H\"{o}lder-type for each of them, uniformly with respect to the parameter $\nu$. From now on, we will do this in the particular case when the parameter $\nu=0$, but the same argument works in general. By Fatou's lemma, we can also restrict the summation in the definition \eqref{3.24} of $\dot{\Pi}^{\delta}_{\mathbb{P}^{k}_{n}}(f_{1},f_{2})$ on arbitrary finite sub-collections $\mathbb{P}$ of $\mathbb{P}^{k}_{n}$, and prove the estimates are unform with respect to different choices of the set $\mathbb{P}$.

\medskip

Therefore, one can reduce the bilinear operator $\dot{\Pi}^{\delta}_{\mathbb{P}^{k}_{n}}$ further to the discrete bilinear model operator $\Pi^{\delta}_{\mathbb{P},n,k}$ defined by
\begin{equation}\label{3.25}
  \Pi^{\delta}_{\mathbb{P},n,k}(f_{1},f_{2})(x):=\sum_{P\in\mathbb{P}\subseteq\mathbb{P}^{k}_{n}}\frac{C^{\delta}_{Q_{P}}}{|I_{P}|^{\frac{1}{2}}}
  \langle f_{1},\Phi^{1}_{P_{1}}\rangle\langle f_{2},\Phi^{2}_{P_{2}}\rangle\Phi^{3}_{P_{3}}(x),
\end{equation}
where $\Phi^{j}_{P_{j}}:=\Phi^{j,0}_{P_{j}}$ for $j=1,2,3$ respectively, $\mathbb{P}\subseteq\mathbb{P}^{k}_{n}$ is arbitrary finite sub-collection of tri-tiles contained in $\mathbb{P}^{k}_{n}$. By \eqref{FC}, one has the following estimates for the Fourier coefficients $C^{\delta}_{Q_{P}}$:
\begin{equation}\label{3.26}
  |C^{\delta}_{Q_{P}}|\lesssim e^{-\delta\sqrt{n}}+n^{-3},
\end{equation}
therefore, we can normalize the model operator $\Pi^{\delta}_{\mathbb{P},n,k}$ by changing the coefficients $C^{\delta}_{Q_{P}}$ into $1$ and only need to prove $L^{p}$ estimates for the normalized model operator with the constants in the bounds depending only on $p$, $p_{1}$, $p_{2}$ and independent of $\delta$, $n$ and $k$.

\medskip

As have discussed above, we now reach a conclusion that the proof of Proposition \ref{single multiplier} can be reduced to proving the following $L^{p}$ estimates for an arbitrary single model operator $\Pi^{n,k}_{\mathbb{P}}$.

\begin{prop}\label{single model}
For every $n=1,2,\cdots$ and $1\leq k\lesssim10^{8}n^{\frac{3}{2}}$, if the finite collection $\mathbb{P}\subseteq\mathbb{P}^{k}_{n}$ is chosen arbitrarily as above, then the discrete model operator $\Pi^{n,k}_{\mathbb{P}}$ defined by
\begin{equation}\label{model operator}
  \Pi^{n,k}_{\mathbb{P}}(f_{1},f_{2})(x):=\sum_{P\in\mathbb{P}\subseteq\mathbb{P}^{k}_{n}}\frac{1}{|I_{P}|^{\frac{1}{2}}}
  \langle f_{1},\Phi^{1}_{P_{1}}\rangle\langle f_{2},\Phi^{2}_{P_{2}}\rangle\Phi^{3}_{P_{3}}(x)
\end{equation}
maps $L^{p_{1}}(\mathbb{R})\times L^{p_{2}}(\mathbb{R})\rightarrow L^{p}(\mathbb{R})$ boundedly for any $1<p_{1},p_{2}\leq\infty$ and $\frac{1}{2}<p<\infty$ satisfying $\frac{1}{p}=\frac{1}{p_{1}}+\frac{1}{p_{2}}$. Moreover, the implicit constants in the bounds depend only on $p_{1}$, $p_{2}$, $p$ and are independent of $n$, $k$ and the particular choice of finite sub-collection $\mathbb{P}$.
\end{prop}

\subsection{Estimates of each single bilinear discrete model operator}

In this section, we prove Proposition \ref{single multiplier} by carrying out the proof of Proposition \ref{single model} for bilinear discrete model operator $\Pi^{n,k}_{\mathbb{P}}$ defined by \eqref{model operator}.

\medskip

We first consider the cases $n\geq2$. It's well known that a standard approach to prove $L^{p}$ estimates for one-parameter $n$-linear operators with singular symbols (e.g., Coifman-Meyer multiplier, $BHT$ and one-parameter paraproducts) is the generic estimates of the corresponding $(n+1)$-linear forms consisting of estimates for different \emph{sizes}, \emph{energies}, or maximal operators $M$ and square operators $S$ (see \cite{Jung,MS,MTT1,MTT2}). To this end, we define the maximal operator $M(f_{1})$ and discretized square operators $S(f_{2})$ and $S(f_{3})$ for functions $f_{1}$, $f_{2}$ and $f_{3}$ by
\begin{equation*}
   M(f_{1})(x):=\sup_{P\in\mathbb{P}\subseteq\mathbb{P}^{k}_{n}}\frac{|\langle f_{1},\Phi^{1}_{P_{1}}\rangle|}{|I_{P}|^{\frac{1}{2}}}\chi_{I_{P}}(x) \,\,\,\,
   \text{and} \,\,\,\, S(f_{i})(x):=\bigg(\sum_{P\in\mathbb{P}\subseteq\mathbb{P}^{k}_{n}}\frac{|\langle f_{i},\Phi^{i}_{P_{i}}\rangle|^{2}}{|I_{P}|}\chi_{I_{P}}(x)\bigg)^{\frac{1}{2}}
\end{equation*}
for $i=2,3$. By the construction of $\mathbf{Q}^{k}_{n}$, one observes that the tri-tiles $P\in\mathbb{P}\subseteq\mathbb{P}^{k}_{n}$ are uniquely determined by its time intervals $I_{P}$, moreover, all the $\{\Phi^{i}_{P_{i}}\}_{P\in\mathbb{P}}$ ($i=1,2,3$) are families of $L^{2}$-normalized bump functions adapted to dyadic intervals $I_{P}$ and $\{\Phi^{i}_{P_{i}}\}_{P\in\mathbb{P}}$ ($i=2,3$) also have the \emph{integral zero} property that $\int_{\mathbb{R}}\Phi^{i}_{P_{i}}(x)dx=0$, thus the maximal operator $M(f_{1})$ and the square operators $S(f_{2})$, $S(f_{3})$ are bounded on every $L^{p}$ space for $1<p<\infty$, and the implicit constants in the bounds will depend only on $p$. The desired $L^{p}$ estimates can be easily deduced from H\"{o}lder estimates in the particular cases $1<p,p_{1},p_{2}<\infty$. Indeed, let $f_{1}\in L^{p_{1}}$, $f_{2}\in L^{p_{2}}$ and $f_{3}\in L^{p'}$ with $\frac{1}{p}=\frac{1}{p_{1}}+\frac{1}{p_{2}}$, $1<p,p_{1},p_{2}<\infty$ and $\|f_{3}\|_{L^{p'}}=1$, then one can derive that
\begin{eqnarray*}
 && \bigg|\int_{\mathbb{R}}\Pi^{n,k}_{\mathbb{P}}(f_{1},f_{2})(x)f_{3}(x)dx\bigg|\lesssim
 \int_{\mathbb{R}}\sum_{P\in\mathbb{P}\subseteq\mathbb{P}^{k}_{n}}
  \frac{|\langle f_{1},\Phi^{1}_{P_{1}}\rangle|}{|I_{P}|^{\frac{1}{2}}}\frac{|\langle f_{2},\Phi^{2}_{P_{2}}\rangle|}{|I_{P}|^{\frac{1}{2}}}
  \frac{|\langle f_{3},\Phi^{3}_{P_{3}}\rangle|}{|I_{P}|^{\frac{1}{2}}}\chi_{I_{P}}(x)dx\\
 \nonumber &\lesssim&\int_{\mathbb{R}}M(f_{1})(x)S(f_{2})(x)S(f_{3})(x)dx
 \lesssim\|M(f_{1})\|_{L^{p_{1}}}\|S(f_{2})\|_{L^{p_{2}}}\|S(f_{3})\|_{L^{p'}}\lesssim\|f_{1}\|_{L^{p_{1}}}\|f_{2}\|_{L^{p_{2}}},
\end{eqnarray*}
and hence
\[\|\Pi^{n,k}_{\mathbb{P}}(f_{1},f_{2})\|_{L^{p}(\mathbb{R})}\lesssim_{p,p_{1},p_{2}}\|f_{1}\|_{L^{p_{1}}(\mathbb{R})}\cdot\|f_{2}\|_{L^{p_{2}}(\mathbb{R})}.\]

\medskip

By the multi-linear interpolations (see \cite{GrTa,Janson,MS,MTT1}) and the symmetry of operators $\Pi^{n,k}_{\mathbb{P}}$, in order to prove $L^{p}$ estimates in the general cases $1<p_{1},p_{2}\leq\infty$ and $\frac{1}{2}<p<\infty$, we only need to prove that the bilinear model operators $\Pi^{n,k}_{\mathbb{P}}$ satisfy the $L^{\frac{1}{2}+\varepsilon,\infty}$ estimates for $\varepsilon>0$ arbitrarily small. To do this, by using the duality lemma for $L^{r,\infty}$ (see \cite{MS,MTT1}) and scaling invariance, we fix $p_{1}$, $p_{2}$ two numbers larger than $1$ and arbitrarily close to $1$, let functions $f_{1}$, $f_{2}$ such that $\|f_{1}\|_{L^{p_{1}}}=\|f_{2}\|_{L^{p_{2}}}=1$ and a measurable set $E\subseteq\mathbb{R}$ satisfying $|E|=1$, our goal is to find a dominative subset $E'\subseteq E$ with comparable measure $|E'|\simeq1$ such that the corresponding trilinear forms $\Lambda^{k}_{\mathbb{P},n}$ satisfy the estimates:
\begin{equation}\label{3.27}
  |\Lambda^{k}_{\mathbb{P},n}(f_{1},f_{2},f_{3})|=\bigg|\sum_{P\in\mathbb{P}\subseteq\mathbb{P}^{k}_{n}}\frac{1}{|I_{P}|^{\frac{1}{2}}}
  \langle f_{1},\Phi^{1}_{P_{1}}\rangle\langle f_{2},\Phi^{2}_{P_{2}}\rangle\langle f_{3},\Phi^{3}_{P_{3}}\rangle\bigg|\lesssim 1,
\end{equation}
where $f_{3}:=\chi_{E'}$.

\medskip

Since one can essentially regard the operators $\Pi^{n,k}_{\mathbb{P}}$ given by \eqref{model operator} as one-parameter paraproducts, we can use the \emph{stopping-time decomposition} arguments based on square operators $S$ and maximal operators $M$ developed by C. Muscalu, J. Pipher, T. Tao and C. Thiele in \cite{MPTT1,MPTT,MPTT2} to prove \eqref{3.27}. We will give a brief proof for \eqref{3.27} here.

\medskip

By using the generic decomposition lemma (Lemma 3.1 in \cite{MPTT2}), one can estimate \eqref{3.27} by
\begin{equation}\label{3.28}
  |\Lambda^{k}_{\mathbb{P},n}(f_{1},f_{2},f_{3})|\lesssim\sum_{\ell\in\mathbb{N}}2^{-1000\ell}\sum_{P\in\mathbb{P}\subseteq\mathbb{P}^{k}_{n}}
  \frac{1}{|I_{P}|^{\frac{1}{2}}}|\langle f_{1},\Phi^{1}_{P_{1}}\rangle||\langle f_{2},\Phi^{2}_{P_{2}}\rangle||\langle f_{3},\Phi^{3,\ell}_{P_{3}}\rangle|,
\end{equation}
where the new bump functions $\Phi^{3,\ell}_{P_{3}}$ also satisfy the \emph{integral zero} property as $\Phi^{3}_{P_{3}}$ but have the additional property that $supp \, (\Phi^{3,\ell}_{P_{3}})\subseteq 2^{\ell}I_{P}$, and hence the square operators $S^{\ell}$ (which are defined in terms of $\Phi^{3,\ell}_{P_{3}}$ instead of $\Phi^{3}_{P_{3}}$) are bounded on any $L^{p}$ ($1<p<\infty$) space as well with the constants in the bounds depending only on $p$. For every $\ell\in\mathbb{N}$, we define the sets as follows:
\begin{equation}\label{3.29}
  \Omega_{-10\ell}:=\{x\in\mathbb{R}: \, M(f_{1})(x)>C2^{10\ell}\}
  \cup\{x\in\mathbb{R}: \, S(f_{2})(x)>C2^{10\ell}\},
\end{equation}
\begin{equation}\label{3.30}
  \widetilde{\Omega}_{-10\ell}:=\{x\in\mathbb{R}: \, M(\chi_{\Omega_{-10\ell}})(x)>\frac{1}{100}\},
\end{equation}
and
\begin{equation}\label{3.31}
  \widetilde{\widetilde{\Omega}}_{-10\ell}:=\{x\in\mathbb{R}: \, M(\chi_{\widetilde{\Omega}_{-10\ell}})(x)>2^{-\ell}\}.
\end{equation}
Finally, we define the exceptional set
\begin{equation}\label{3.32}
  U:=\bigcup_{\ell\in\mathbb{N}}\widetilde{\widetilde{\Omega}}_{-10\ell}.
\end{equation}
It is clear that $|U|<\frac{1}{10}$ if $C$ is a large enough constant, which we fix from now on. Then, we define the dominative subset $E':=E\setminus U$ and observe that $|E'|\simeq1$.

\medskip

Now we fix $\ell\in\mathbb{N}$ and consider the corresponding inner sum in \eqref{3.28}. One easily observes that it's enough for us to consider the sub-collection $\mathbb{P}_{\ell}\subseteq\mathbb{P}\subseteq\mathbb{P}^{k}_{n}$ consisting of all tri-tiles $P\in\mathbb{P}$ with time intervals $I_{P}$ satisfying $I_{P}\bigcap\widetilde{\Omega}^{c}_{-10\ell}\neq\emptyset$ in the inner sum of \eqref{3.28} (other parts are equal to $0$), it follows that $|I_{P}\bigcap\Omega_{-10\ell}|\leq\frac{1}{100}|I_{P}|$ for such tri-tiles $P\in\mathbb{P}_{\ell}$.

\medskip

Now we define three different decomposition procedures for functions $f_{1}$, $f_{2}$ and $f_{3}$ respectively. First, define
\[\Omega^{1}_{n_{1}}:=\{x\in\mathbb{R}|M(f_{1})(x)>C2^{-n_{1}}\}, \,\,\,\,\,\, \Omega^{2}_{n_{2}}:=\{x\in\mathbb{R}|S(f_{2})(x)>C2^{-n_{2}}\}\]
and
$$\mathbb{P}^{i}_{n_{i},\ell}:=\bigg\{P\in\mathbb{P}_{\ell}\bigg||I_{P}\bigcap\Omega^{i}_{n_{i}}|>\frac{|I_{P}|}{100}, \, |I_{P}\bigcap\Omega^{i}_{n_{i}-1}|\leq\frac{|I_{P}|}{100}\bigg\}$$
for $i=1,2$ and every $n_{1},n_{2}\geq-10\ell$, which produces the sets $\{\Omega^{1}_{n_{1}}\}_{n_{1}\geq-10\ell}$, $\{\Omega^{2}_{n_{2}}\}_{n_{2}\geq-10\ell}$, $\{\mathbb{P}^{1}_{n_{1},\ell}\}_{n_{1}\geq-10\ell}$ and $\{\mathbb{P}^{2}_{n_{2},\ell}\}_{n_{2}\geq-10\ell}$ with $\mathbb{P}^{1}_{-10\ell,\ell}=\mathbb{P}^{2}_{-10\ell,\ell}=\emptyset$, and one has $\mathbb{P}_{\ell}=\bigcup_{n_{i}\geq-10\ell}\mathbb{P}^{i}_{n_{i},\ell}$ for $i=1,2$. Then, we choose $N>0$ large enough such that for every $P\in\mathbb{P}$, one has $|I_{P}\cap\Omega^{3}_{-N}|\leq\frac{|I_{P}|}{100}$, where $\Omega^{3}_{-N}:=\{x\in\mathbb{R}|S^{\ell}(f_{3})(x)>C2^{N}\}$. For every $n_{3}\geq-N$, we define sets
\[\Omega^{3}_{n_{3}}:=\{x\in\mathbb{R}|S^{\ell}(f_{3})(x)>C2^{-n_{3}}\}\]
and
$$\mathbb{P}^{3}_{n_{3},\ell}:=\bigg\{P\in\mathbb{P}_{\ell}\bigg||I_{P}\bigcap\Omega^{3}_{n_{3}}|>\frac{|I_{P}|}{100}, \, |I_{P}\bigcap\Omega^{3}_{n_{3}-1}|\leq\frac{|I_{P}|}{100}\bigg\},$$
which produces the sets $\{\Omega^{3}_{n_{3}}\}_{n_{3}\geq-N}$ and $\{\mathbb{P}^{3}_{n_{3},\ell}\}_{n_{3}\geq-N}$ with $\mathbb{P}^{3}_{-N,\ell}=\emptyset$, and one has $\mathbb{P}_{\ell}=\bigcup_{n_{3}\geq-N}\mathbb{P}^{3}_{n_{3},\ell}$.

\medskip

Now we define sub-collection $\mathbb{P}^{n_{1},n_{2},n_{3}}_{\ell}:=\mathbb{P}^{1}_{n_{1},\ell}\cap\mathbb{P}^{2}_{n_{2},\ell}\cap\mathbb{P}^{3}_{n_{3},\ell}$ and $\Omega^{\ell}_{n_{1},n_{2},n_{3}}:=\bigcup_{P\in\mathbb{P}^{n_{1},n_{2},n_{3}}_{\ell}}I_{P}$. The inner sum in the right-hand side of \eqref{3.28} can be estimated by
\begin{eqnarray}\label{3.33}
&& \sum_{\substack{n_{1},n_{2}>-10\ell, \\ n_{3}>-N}}\sum_{P\in\mathbb{P}^{n_{1},n_{2},n_{3}}_{\ell}}
\frac{|\langle f_{1},\Phi^{1}_{P_{1}}\rangle|}{|I_{P}|^{\frac{1}{2}}}\cdot\frac{|\langle f_{2},\Phi^{2}_{P_{2}}\rangle|}{|I_{P}|^{\frac{1}{2}}}
\cdot\frac{|\langle f_{3},\Phi^{3,\ell}_{P_{3}}\rangle|}{|I_{P}|^{\frac{1}{2}}}\cdot\bigg|I_{P}\cap\big(\bigcup_{i=1}^{3}\Omega^{i}_{n_{i}-1}\big)^{c}\bigg|\\
\nonumber &\lesssim& \sum_{\substack{n_{1},n_{2}>-10\ell, \\ n_{3}>-N}}\int_{(\bigcup_{i=1}^{3}\Omega^{i}_{n_{i}-1})^{c}\cap\Omega^{\ell}_{n_{1},n_{2},n_{3}}}
M(f_{1})(x)S(f_{2})(x)S^{\ell}(f_{3})(x)dx\\
\nonumber &\lesssim& \sum_{\substack{n_{1},n_{2}>-10\ell, \\ n_{3}>-N}}2^{-(n_{1}+n_{2}+n_{3})}|\Omega^{\ell}_{n_{1},n_{2},n_{3}}|.
\end{eqnarray}
By the boundedness of the operators $M$, $S$ and $S^{\ell}$, one easily get the estimates
\begin{equation}\label{3.34}
  |\Omega^{\ell}_{n_{1},n_{2},n_{3}}|\lesssim 2^{n_{1}p_{1}}, \, 2^{n_{2}p_{2}}, \, 2^{n_{3}\mu}
\end{equation}
for any $\mu\geq1$, as a consequence, combining this with the estimates \eqref{3.28} and \eqref{3.33}, we derive the following estimates for the trilinear form:
\begin{eqnarray}\label{3.35}
|\Lambda^{k}_{\mathbb{P},n}(f_{1},f_{2},f_{3})|&\lesssim&\sum_{\ell\in\mathbb{N}}2^{-1000\ell}
\{\sum_{\substack{n_{1},n_{2}>-10\ell, \\ n_{3}>0}}2^{-(1-p_{1}\theta_{1})n_{1}}2^{-(1-p_{2}\theta_{2})n_{2}}2^{-(1-\mu\theta_{3})n_{3}}\\
\nonumber &&\quad\quad\quad\quad\quad +\sum_{\substack{n_{1},n_{2}>-10\ell, \\ -N<n_{3}<0}}2^{-(1-p_{1}\theta'_{1})n_{1}}2^{-(1-p_{2}\theta'_{2})n_{2}}2^{-(1-\mu\theta'_{3})n_{3}}\}\\
\nonumber &\lesssim& \sum_{\ell\in\mathbb{N}}2^{-100\ell}\lesssim 1,
\end{eqnarray}
where the parameter $\theta_{i}$ and $\theta'_{i}$ are chosen appropriately, which implies the desired estimate \eqref{3.27}.

\medskip

As to the cases $n=1$ and $1\leq k\lesssim10^{8}$, we will perform the discretized square operators $S$ on the functions $f_{1}$ and $f_{2}$ and the maximal operator $M$ on function $f_{3}$, and estimates the corresponding trilinear form $\Lambda^{k}_{\mathbb{P},n}(f_{1},f_{2},f_{3})$ by $S(f_{1})$, $S(f_{2})$ and $M(f_{3})$, the rest of the proof are completely similar to the $n\geq2$ cases.

This completes the proof of Proposition \ref{single model}, and also concludes the proof of Proposition \ref{single multiplier} at the same time.

\subsection{Conclusions}

We can deduce from Proposition \ref{single multiplier} and \eqref{3.13} the following $L^{p}$ estimates for the bilinear operator $T_{m_{\mathcal{A}}}$ (associated with symbol $m$ smoothly truncated on the cone $\mathcal{A}$):
\begin{eqnarray}\label{conclusion}
 && \|T_{m_{\mathcal{A}}}(f_{1},f_{2})\|_{L^{p}(\mathbb{R})}\lesssim\sum_{n=1}^{\infty}\sum_{k=1}^{\sim10^{8}n^{\frac{3}{2}}}
  \|T^{k}_{m_{\mathcal{A}},n}(f_{1},f_{2})\|_{L^{p}(\mathbb{R})}\\
 \nonumber &\lesssim_{p,p_{1},p_{2}}&\sum_{n=1}^{\infty}\sum_{k=1}^{\sim10^{8}n^{\frac{3}{2}}}\{e^{-\delta\sqrt{n}}+n^{-3}\}
  \|f_{1}\|_{L^{p_{1}}(\mathbb{R})}\|f_{2}\|_{L^{p_{2}}(\mathbb{R})}
  \lesssim_{p,p_{1},p_{2},\delta}\|f_{1}\|_{L^{p_{1}}(\mathbb{R})}\cdot\|f_{2}\|_{L^{p_{2}}(\mathbb{R})}
\end{eqnarray}
for any $\frac{1}{p}=\frac{1}{p_{1}}+\frac{1}{p_{2}}$ with $1<p_{1}, \, p_{2}\leq\infty$ and $\frac{1}{2}<p<\infty$, where the implicit constants in the bounds depend only on $p_{1}$, $p_{2}$, $p$, $\delta$ and tend to infinity as $\delta\rightarrow0$. One can observe that \eqref{conclusion} also implies the $L^{p}$ estimates of the original bilinear operator $T_{m}$ for $p>\frac{1}{2}$.

This concludes the proof of Theorem \ref{main1}.\\

\end{document}